\documentclass[11pt]{article}
\usepackage{amsmath,amssymb,amstext,dsfont,fancyvrb,float,fontenc,graphicx,subfigure, theorem}

\title{\Large\bf On the Riemann-Hilbert Problem IV}
\author{\sc Vladimir Ryazanov}
\date{}

\cleardoublepage 
\usepackage[strict]{changepage}
\def\abstractname{Abstract -}   
\def\abstract{\begin{adjustwidth}{1cm}{1cm} \par    \footnotesize \noindent {\bf \abstractname}
\def\endabstract{ \end{adjustwidth} \smallskip }}

{\theorembodyfont{\itshape}\newtheorem{theorem}{Theorem}[section]}
{\theorembodyfont{\itshape}\newtheorem{proposition}{Proposition}[section]}
{\theorembodyfont{\itshape}}
{\theorembodyfont{\itshape}}
{\theorembodyfont{\itshape}\newtheorem{corollary}{Corollary}[section]}
{\theorembodyfont{\rm}}
{\theorembodyfont{\rm}\newtheorem{remark}{Remark}[section]}
{\theorembodyfont{\rm}}
{\theorembodyfont{\rm }}

\begin{document}
\maketitle
\vskip 1.5em


\vskip 1.5em
 \begin{abstract}
It is proved the existence of solutions for the Riemann-Hilbert
problem in the fairly general settings of arbitrary Jordan domains,
mea\-su\-rab\-le coefficients and measurable boundary data. The
theorem is formulated in terms of harmonic measure and principal
asymptotic values. It is also given the corresponding reinforced
criterion for domains with arbitrary rectifiable boundaries stated
in terms of the natural parameter and non\-tan\-gen\-tial limits.
Furthemore, it is shown that the dimension of the spaces of
solutions is infinite.
 \end{abstract}

\begin{keywords} Riemann-Hilbert problem, Jordan domains, harmonic measure, principal
asymptotic values, nontangential limits.
\end{keywords}

\begin{MSC}
primary   31A05, 31A20, 31A25, 31B25, 35Q15; se\-con\-da\-ry 30E25,
31C05, 34M50, 35F45.
\end{MSC}

\section{Introduction}

Boundary value problems for analytic functions are due to the well-known Riemann dissertation (1851) contained a
general setting of a problem on finding analytic functions with a connection between its real and imaginary
parts on the boundary. However, it has contained no concrete boundary value problems.

\medskip

The first concrete problem of such a type has been proposed by Hilbert (1904) and called by the Hilbert problem
or the Riemann-Hilbert problem. That consists in finding an analytic function $f$ in a domain bounded by a
rectifiable Jordan curve $C$ with the linear boundary condition
\begin{equation}\label{eqLIN}
\lim\limits_{z\to\zeta}\ \mathrm {Re}\
\{\overline{\lambda(\zeta)}\cdot f(z)\}\ =\ \varphi(\zeta)
\quad\quad\quad\ \ \ \forall \ \zeta\in C
\end{equation}
where it was assumed by him that the functions $\lambda$ and $\varphi$ are continuously differentiable with
respect to the natural parameter $s$ on $C$ and, moreover, $|\lambda|\ne 0$ everywhere on $C$. Hence without
loss of generality one can assume that $|\lambda|\equiv 1$ on $C$.

\medskip

The first way for solving this problem based on the theory of singular integral equations was given by Hilbert
(1904), see \cite{H1}. This attempt was not quite successful because of the theory of singular integral
equations has been not yet enough developed at that time. However, just that way became the main approach in
this research direction with important contributions of Georgian and Russian mathematicians and mechanicians,
see e.g. \cite{G}, \cite{M} and \cite{V}. In particular, the existence of solutions to this problem was  in that
way proved for H\"older continuous $\lambda$ and $\varphi$. But subsequent weakening conditions on $\lambda$ and
$\varphi$ led to strengthening conditions on the contour $C$, say to the Lyapunov curves or the Radon condition
of bounded rotation or even to smooth curves.

\medskip

However, Hilbert (1905) has proposed the second way for solving his problem in setting to (\ref{eqLIN}) above
based on the reduction it to solving the corresponding two Dirichlet problems, see e.g. \cite{H2}. The goal of
this paper is to show that this approach is more simple and leads to perfectly general results. That requests to
apply some fundamental concepts and facts related to the Dirichlet problem.

\bigskip


\section{The case of the unit circle}

The following brilliant result of Frederick Gehring is key for our goals, see \cite{Ge}.

\medskip
\begin{proposition}\label{prGe}
{\it\, Let $\varphi(\vartheta)$ be real, measurable, almost everywhere finite and have the period $2\pi$. Then
there exists a function $u(z)$, harmonic in $|z|<1$, such that $u(z)\to \varphi(\vartheta)$ for a.e. $\vartheta$
as $z\to e^{i\vartheta}$ along any nontangential path.}\end{proposition}

\medskip

Since the Gehring proof is very short and nice and has a common interest, we give it for completeness here.

\medskip

\begin{proof}
By a theorem of Lusin, see e.g. Theorem VII(2.3) in \cite{S}, p. 217, we can find a continuous function
$\Phi(\vartheta)$ such that $\Phi^{\prime}(\vartheta)=\varphi(\vartheta)$ for a.e. $\vartheta$. Let
$$
U(re^{i\vartheta})\ =\ \frac{1}{2\pi}\ \int\limits_{0}\limits^{2\pi}\frac{1-r^2}{1-2r\cos(\vartheta-t)+r^2}\
\Phi(t)\ dt
$$
for $r<1$. Next, by the well-known result due to Fatou, see e.g.
3.441 in \cite{Z}, p. 53, $\frac{\partial}{\partial\vartheta}\
U(z)\to \Phi^{\prime}(\vartheta)$ as $z\to e^{i\vartheta}$ along any
nontangential path whenever $\Phi^{\prime}(\vartheta)$ exists. Thus,
the conclusion follows for the function $u(z)\ =\
\frac{\partial}{\partial\vartheta}\ U(z)$. \end{proof}

\medskip

\begin{remark}\label{remOS2.1.a} Recall also the preceding result of W.
Kaplan on the existence of a harmonic function $u(z)$ with the radial limits $\varphi(\vartheta)$ a.e., see
\cite{K}.
\end{remark}

It is known that every harmonic function $u(z)$ in $\mathbb D = \{ z\in\mathbb C : |z|<1\}$ has a conjugate
function $v(z)$ such that $f(z)=u(z)+iv(z)$ is an analytic function in $\mathbb D$. Hence we have the following
consequence of Proposition \ref{prGe}.

\medskip
\begin{corollary}\label{corGe}
{\it\, Under the conditions of Proposition \ref{prGe}, there exists an analytic function $f$ in $\mathbb D$ such
that $\mathrm Re\ f(z)\to \varphi(\vartheta)$ for a.e. $\vartheta$ as $z\to e^{i\vartheta}$ along any
nontangential path.}
\end{corollary}

Note that the boundary values of the conjugate function $v$ cannot be prescribed arbitrarily and simultaneously
with the boundary values of $u$ because $v$ is uniquely determined by $u$ up to an additive constant.

\medskip

Denote by $h^p$, $p\in(0,\infty)$, the class of all harmonic functions $u$ in $\mathbb D$ with
$$
\sup\limits_{r\in(0,1)}\ \left\{\int\limits_{0}\limits^{2\pi}\ |u(re^{i\vartheta})|^p\
d\vartheta\right\}^{\frac{1}{p}}\ <\ \infty\ .
$$
It is clear that $h^p\subseteq h^{p^{\prime}}$ for all $p>p^{\prime}$ and, in particular, $h^p\subseteq h^{1}$
for all $p>1$. It is important that every function in the class $h^1$ has a.e. nontangential boundary limits,
see e.g. Corollary IX.2.2 in \cite{Go}.

\medskip

Note also that $v\in h^p$ whenever $u\in h^p$ for all $p>1$ by the M. Riesz theorem, see \cite{R}. Generally
speaking, this fact is not trivial but it follows immediately for $p=2$ from the Parseval equality. The latter
will be sufficient for our goals.

\medskip
\begin{theorem}{}\label{thRHD}{\it\, Let $\lambda:\partial\mathbb D\to\mathbb C$,
$|\lambda (\zeta)|\equiv 1$, and $\varphi:\partial\mathbb
D\to\mathbb R$ be measurable functions. Then there exist analytic
functions $f:\mathbb D\to\mathbb C$ such that along any
nontangential path
\begin{equation}\label{eqLIMD} \lim\limits_{z\to\zeta}\ \mathrm
{Re}\ \{\overline{\lambda(\zeta)}\cdot f(z)\}\ =\ \varphi(\zeta)
\quad\quad\quad\mbox{for}\ \ \mbox{a.e.}\ \ \
\zeta\in\partial\mathbb D\ .\end{equation}} \end{theorem}

\medskip

\begin{proof}
First, consider the function $\alpha(\zeta)=\arg\lambda(\zeta)$
where $\arg \omega$  is the principal value of the argument of
$\omega\in\mathbb C$ with $|\omega|=1$, i.e., the unique number
$\alpha\in (-\pi,\pi]$ such that $\omega=e^{i\alpha}$. Note that the
function $\arg \omega$ is continuous on $\partial\mathbb
D\setminus\{{-1}\}$ and the sets $\lambda^{-1}(\partial\mathbb
D\setminus\{{-1}\})$ and $\lambda^{-1}({-1})$ are measurable because
the function $\lambda(\zeta)$ is measurable. Thus, the function
$\alpha(\zeta)$ is measurable on $\partial\mathbb D$. Furthermore,
$\alpha\in L^{\infty}(\partial\mathbb D)$ because
$|\alpha(\zeta)|\le \pi$ for all $\zeta\in\partial\mathbb D$. Hence
\begin{equation}\label{eqPOISSON} g(z)\ \colon =\ \frac{1}{2\pi i}\ \int\limits_{\partial\mathbb
D}\alpha(\zeta)\ \frac{z+\zeta}{z-\zeta}\
 \frac{d\zeta}{\zeta}\ , \ \ \ \ \ z\in\mathbb D\ ,
\end{equation} is an analytic function in $\mathbb D$ with
$u(z)={\mathrm Re}\ g(z)\to\alpha(\zeta)$ as $z\to\zeta$ along any
nontangential path in $\mathbb D$ for a.e. $\zeta\in\partial\mathbb
D$, see, e.g., Corollary IX.1.1 in \cite{Go} and  Theorem I.E.1 in
\cite{Ko}. Denote ${\cal A}(z)=\exp\{ig(z)\}$ that is an analytic
function.

Since $\alpha\in L^{\infty}(\partial\mathbb D)$, we have that $u\in
h^p$ for all $p\ge 1$, see, e.g., Theorem IX.2.3 in \cite{Go}, and
then $v={\mathrm Im}\ g \in h^p$ for all $p\ge 1$  by the theorem of
M. Riesz. Hence there exists a function $\beta:\partial\mathbb
D\to\mathbb R$, $\beta\in L^p$, for all $p\ge 1$ such that
$v(z)\to\beta(\zeta)$ as $z\to\zeta$ for a.e.
$\zeta\in\partial\mathbb D$ along any nontangential path, see e.g.
Theorem IX.2.3 and Corollary IX.2.2 in \cite{Go}. Thus, by Corollary
\ref{corGe} there exists an analytic function ${\cal B}:\mathbb
D\to\mathbb C$ such that $\mathrm {Re}\ {\cal B}(z) = B(\zeta)\
\colon =\varphi(\zeta)\cdot \exp\{{\beta(\zeta)}\}$ as $z\to\zeta$
along any nontangential path for a.e. $\zeta\in\partial\mathbb D$.
Finally, elementary calculations show that one of the desired
analytic functions in (\ref{eqLIMD}) is $f={\cal A}\cdot{\cal B}$.
\end{proof}

\begin{remark}\label{remOS2.1.A} As it follows from the formula
(\ref{eqPOISSON}), the first analytic function ${\cal A}$ in the
proof is calculated in the explicit form. The function
$\beta:\partial\mathbb D\to\mathbb R$ in the proof can also
explicitly be calculated by the following formula, see, e.g.,
Theorem I.E.4.1 in \cite{Ko}, for a.e. $\zeta\in\partial\mathbb D$
\begin{equation}\label{eqHILBERT} \beta(\zeta)\ \colon =\
\lim\limits_{\varepsilon\to +0}\ \frac{1}{\pi}\
\int\limits_{\varepsilon}\limits^{\pi} \frac{\alpha(\zeta
e^{-it})-\alpha(\zeta e^{it})}{2\ {\rm tg}\ \frac{t}{2}}\ dt\ .
\end{equation}
The second analytic function ${\cal B}$ in the proof is equal to
$\frac{\partial}{\partial\vartheta}\ G(z)$, $z=re^{i\vartheta}$,
with
\begin{equation}\label{eqLUSIN} G(z)\ \colon =\ \frac{1}{2\pi i}\ \int\limits_{\partial\mathbb
D}\Phi(\zeta)\ \frac{z+\zeta}{z-\zeta}\
 \frac{d\zeta}{\zeta}\ , \ \ \ \ \ z\in\mathbb D\ ,
\end{equation}
where $\Phi:\partial\mathbb D\to\mathbb R$ is a continuous function
such that $\frac{\partial}{\partial\vartheta}\
\Phi(\zeta)=B(\zeta)$, $\zeta=e^{i\vartheta}$, for a.e.
$\vartheta\in[0,2\pi]$, see the nontrivial construction of Theorem
VII(2.3) in \cite{S}.
\end{remark}

\bigskip

\section{The case of a rectifiable Jordan curve}

\begin{theorem}{}\label{thRHR}{\it\, Let $D$ be a Jordan domain in
$\mathbb C$ with a rectifiable boundary and let $\lambda:\partial
D\to\mathbb C$, $|\lambda (\zeta)|\equiv 1$ and $\varphi:\partial
D\to\mathbb R$ be measurable functions with respect to the natural
parameter on $\partial D$. Then there exist analytic functions
$f:\mathbb D\to\mathbb C$ such that along any nontangential path
\begin{equation}\label{eqLIM} \lim\limits_{z\to\zeta}\ \mathrm {Re}\
\{\overline{\lambda(\zeta)}\cdot f(z)\}\ =\ \varphi(\zeta)
\quad\quad\quad\mbox{for}\ \ \mbox{a.e.}\ \ \ \zeta\in\partial D
\end{equation}
with respect to the natural parameter on $\partial D$.}
\end{theorem}

\medskip

\begin{proof} This case is reduced to the case of the unit disk $\mathbb D$ in the following
way. First, by the Riemann theorem, see e.g. Theorem II.2.1 in \cite{Go}, there exists a conformal mapping
$\omega$ of any Jordan domain $D$ onto $\mathbb D$. By the Caratheodory (1912) theorem $\omega$ can be extended
to a homeomorphisms of $\overline D$ onto $\overline{\mathbb D}$ and, if $\partial D$ is rectifiable, then by
the theorem of F. and M. Riesz (1916) $\mathrm{length}\ \omega^{-1}(E)=0$ whenever $E\subset\partial\mathbb D$
with $|E|=0$, see e.g. Theorem II.C.1 and Theorems II.D.2 in \cite{Ko}. Conversely, by the Lavrentiev (1936)
theorem $|\omega({\cal E})|=0$ whenever ${\cal E}\subset\partial D$ and $\mathrm{length}\ {\cal E}=0$, see
\cite{L}, see also the point III.1.5 in \cite{P}.

\medskip

Hence $\omega$ and $\omega^{-1}$ transform measurable sets into measurable sets. Indeed, every measurable set is
the union of a sigma-compact set and a set of measure zero, see e.g. Theorem III(6.6) in \cite{S}, and
continuous mappings transform compact sets into compact sets. Thus, a function $\varphi:\partial D\to\mathbb R$
is measurable with respect to the natural parameter on $\partial D$ if and only if the function
$\Phi=\varphi\circ\omega^{-1}:\partial\mathbb D\to\mathbb R$ is measurable with respect to the linear measure on
$\partial\mathbb D$.

\medskip

By the Lindel\"of (1917) theorem, see e.g. Theorem II.C.2 in \cite{Ko}, if $\partial D$ has a tangent at a point
$\zeta$, then $\arg\ [\omega(\zeta)-\omega(z)]-\arg\ [\zeta-z]\to\mathrm const$ as $z\to\zeta$. In other words,
the conformal images of sectors in $D$ with a vertex at $\zeta$ is asymptotically the same as sectors in
$\mathbb D$ with a vertex at $w=\omega(\zeta)$. Thus, nontangential paths in $D$ are transformed under $\omega$
into nontangential paths in $\mathbb D$. Finally, a rectifiable Jordan curve has a tangent a.e. with respect to
the natural parameter and, thus, Theorem \ref{thRHR} follows from Theorem \ref{thRHD}. \end{proof}\ 

\medskip

In particular, choosing $\lambda\equiv 1$ in (\ref{eqLIM}), we obtain the following statement.

\medskip
\begin{proposition}\label{prDR}
{\it\, Let $D$ be a domain in $\mathbb C$ bounded by a rectifiable Jordan curve and $\varphi:\partial
D\to\mathbb R$ be measurable. Then there exists an analytic function $f:D\to\mathbb C$ such that
\begin{equation}\label{eqLIMDRA} \lim\limits_{z\to\zeta}\ \mathrm
{Re}\
 f(z)\ =\ \varphi(\zeta)
\quad\quad\quad\mbox{for}\ \ \mbox{a.e.}\ \ \ \zeta\in\partial D
\end{equation}
with respect to the natural parameter on $\partial D$ along any nontangential path.}\end{proposition}

\begin{corollary}\label{corDR}
{\it\, Under the conditions of Proposition \ref{prDR}, there exists a harmonic function $u$ in $D$ such that
$\mathrm u(z)\to \varphi(\zeta)$ for a.e. $\zeta\in\partial D$ as $z\to\zeta$ along any nontangential path.}
\end{corollary}

\bigskip

\section{The case of an arbitrary Jordan curve}

The conceptions of a harmonic measure introduced by R. Nevanlinna in
\cite{N} and a principal asymptotic value based on one nice result
of F. Bagemihl \cite{B} make possible with a great simplicity and
generality to formulate the existence theorems for the Dirichlet and
Riemann-Hilbert problems.

\medskip

First of all, given a measurable set $E\subseteq\partial\mathbb D$
and a point $z\in\mathbb D$, a {\it harmonic measure} of $E$ at $z$
relative to $\mathbb D$ is the value at $z$ of the harmonic function
$u$ in $\mathbb D$ with the boundary values $1$ a.e. on $E$ and $0$
a.e on $\partial\mathbb D\setminus E$, see Proposition \ref{prGe}.
In particular, by the mean value theorem for harmonic functions, the
harmonic measure of $E$ at $0$ relative to $\mathbb D$ is equal to
$|E|/2\pi$. In general, the geometric sense of the harmonic measure
of $E$ at $z_0$ relative to $\mathbb D$ is the angular measure of
view of $E$ from the point $z_0$ in radians divided by $2\pi$.

\medskip

Since the harmonic measure zero is invariant under conformal
mappings between Jordan domains, given a Jordan domain $D$, a set
${\cal E}\subseteq \partial D$ will be called measurable with
respect to harmonic measures in $D$ if $E=\omega({\cal E})$ is
measurable with respect to the linear measure on $\partial\mathbb D$
where $\omega$ is a conformal mapping of $D$ onto the unit disk
$\mathbb D$, cf. the proof of Theorem \ref{thRHR}. Correspondingly,
the harmonic measure of ${\cal E}\subseteq
\partial D$ at $z_0\in D$ relative to $D$ is the harmonic measure of
$\omega_0({\cal E})$ at $0$ relative to $\mathbb D$ where $\omega_0$
is a conformal mapping of $D$ onto $\mathbb D$ with the
normalization $\omega_0(z_0)=0$, i.e., the quantity $|\omega_0({\cal
E})|/2\pi$.

\begin{figure}[h]
\centerline{\includegraphics[scale=1.0]{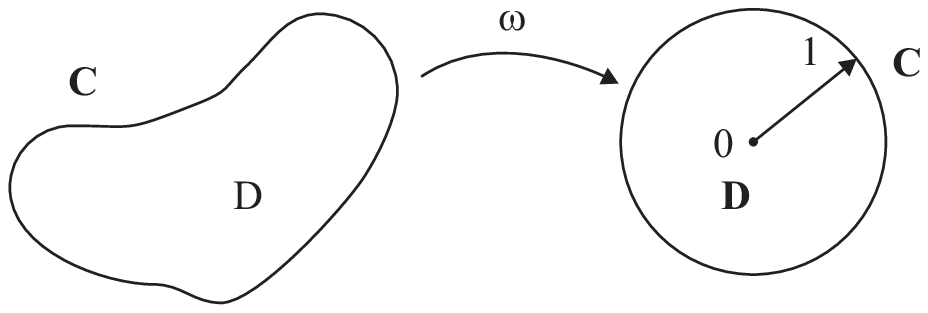}}%
\centerline{Figure 1. The Riemann mapping}
\end{figure}

\medskip

Next, a Jordan curve generally speaking has no tangents. Hence we
need a replacement for the notion of a nontangential limit. In this
connection, recall Theorem 2 in \cite{B}, see also Theorem III.1.8
in \cite{No}, stating that, for any function $\Omega:\mathbb
D\to\overline{\mathbb C}$, for all pairs of arcs $\gamma_1$ and
$\gamma_2$ in $\Bbb D$ terminating at $\zeta\in\partial\mathbb D$,
except a countable set of $\zeta\in\partial\mathbb D$,
\begin{equation}\label{eqBA}
C(\Omega,\gamma_1)\ \cap\ C(\Omega,\gamma_2)\ \neq\ \varnothing
\end{equation}
 where $C(\Omega,\gamma)$ denotes the
{\it cluster set of $\Omega$ at $\zeta$ along $\gamma$}, i.e.,
$$
C(\Omega,\gamma)\ =\ \{ w\in\overline{\mathbb C}\ :\ \Omega(z_n)\to
w,\ z_n\to\zeta ,\ z_n\in\gamma\}\ .
$$
Immediately by the theorems of Riemann and Caratheodory, this result
is extended to an arbitrary Jordan domain $D$ in $\mathbb C$. Given
a function $\Omega: D\to\overline{\mathbb C}$ and $\zeta\in\partial
D$, denote by $P(\Omega , \zeta)$ the intersection of all cluster
sets $C(\Omega,\gamma)$ for arcs $\gamma$ in $D$ terminating at
$\zeta$. Later on, we call the points of the set $P(\Omega , \zeta)$
{\it principal asymptotic values} of $\Omega$ at $\zeta$. Note that,
if $\Omega$ has a limit along at least one arc in $D$ terminating at
a point $\zeta\in\partial D$ with the prpoperty (\ref{eqBA}), then
the principal asymptotic value is unique.

\medskip

Thus, by the Bagemihl theorem, we obtain the following result
directly from Theorem \ref{thRHD}.

\medskip
\begin{theorem}{}\label{thRH1}{\it\, Let $D$ be a Jordan domain in
$\mathbb C$ and let $\lambda:\partial D\to\mathbb C$, $|\lambda
(\zeta)|\equiv 1$, and $\varphi:\partial D\to\mathbb R$ be
measurable functions with respect to harmonic measures in $D$. Then
there exist analytic functions $f:\mathbb D\to\mathbb C$ such that
\begin{equation}\label{eqLIMA} \lim\limits_{z\to\zeta}\ \mathrm {Re}\
\{\overline{\lambda(\zeta)}\cdot f(z)\}\ =\ \varphi(\zeta)
\quad\quad\quad\mbox{for}\ \ \mbox{a.e.}\ \ \ \zeta\in\partial D
\end{equation}
with respect to harmonic measures in $D$ in the sense of the unique
principal asymptotic value.}
\end{theorem}

\medskip

In particular, choosing $\lambda\equiv 1$ in (\ref{eqLIMA}), we
obtain the following consequence.

\medskip
\begin{proposition}\label{prDA}
{\it\, Let $D$ be a Jordan domain and $\varphi:\partial D\to\mathbb
R$ be measurable with respect to harmonic measures in $D$. Then
there exists an analytic function $f:D\to\mathbb C$ such that
\begin{equation}\label{eqLIMDA} \lim\limits_{z\to\zeta}\ \mathrm
{Re}\
 f(z)\ =\ \varphi(\zeta)
\quad\quad\quad\mbox{for}\ \ \mbox{a.e.}\ \ \ \zeta\in\partial D
\end{equation}
with respect to harmonic measures in $D$ in the sense of the unique
principal asymptotic value.}\end{proposition}

\medskip
\begin{corollary}\label{corDA}
{\it\, Under the conditions of Proposition \ref{prDA}, there exists
a harmonic function $u$ in $D$ such that in the same sense
\begin{equation}\label{eqLIMDAH} \lim\limits_{z\to\zeta}\ \mathrm
u(z)\ =\ \varphi(\zeta) \quad\quad\quad\mbox{for}\ \ \mbox{a.e.}\ \
\ \zeta\in\partial D\ .
\end{equation}}\end{corollary}

\begin{remark}\label{remOS2.1.z}
In view of the theorems of Riemann and Caratheodory, this approach
makes possible also to formulate the corresponding theorems for
ar\-bi\-tra\-ry simply connected domains $D$ in $\mathbb C$ having
at least 2 boundary points. The only difference is that the
functions $\lambda$ and $\varphi$ should be given as functions of
prime ends of $D$ but not of points of $\partial D$ and harmonic
measures of sets of prime ends are given through the natural
one-to-one correspondence between the prime ends of $D$ and the
boundary points of $\mathbb D$ under Riemann mappings $\omega :
D\to\mathbb D$, see e.g. \cite{CL}.
\end{remark}

\bigskip

\section{On the dimension of spaces of solutions}

By the Lindel\"of maximum principle, see e.g. Lemma 1.1 in
\cite{GM}, it follows the uniqueness theorem for the Dirichlet
problem in the class of bounded harmonic functions on the unit disk
$\mathbb D = \{ z\in\mathbb{C}: |z|<1\}$. In general there is no
uniqueness theorem in the Dirichlet problem for the Laplace
equation.

\medskip

Furthermore, it was proved in \cite{R1} that the space of all
harmonic functions in $\mathbb D$ with nontangential limit $0$ at
a.e. point of $\partial\mathbb D$ has the infinite dimension.
Thanking to it, the statements on the infinite dimension of the
space of solutions can be extended to the Riemann-Hilbert problem
because above we have reduced the latter to the corresponding two
Dirichlet problems.

\medskip

\begin{theorem}{}\label{thDIM1}{\it\, Let $\lambda:\partial\mathbb D\to\mathbb C$,
$|\lambda (\zeta)|\equiv 1$, and $\varphi:\partial\mathbb
D\to\mathbb R$ be measurable functions. Then the space of all
analytic functions $f:\mathbb D\to\mathbb C$ such that along any
nontangential path
\begin{equation}\label{eqLIMDIM} \lim\limits_{z\to\zeta}\ \mathrm
{Re}\ \{\overline{\lambda(\zeta)}\cdot f(z)\}\ =\ \varphi(\zeta)
\quad\quad\quad\mbox{for}\ \ \mbox{a.e.}\ \ \
\zeta\in\partial\mathbb D \end{equation}} has the infinite
dimension.
\end{theorem}

\medskip

\begin{proof} Let $u:\mathbb D\to\mathbb R$ be a harmonic function with nontangential limit $0$ at
a.e. point of $\partial\mathbb D$ from Theorem 1 in \cite{R1}. Then
there is the unique harmonic function $v:\mathbb D\to\mathbb R$ with
$v(0)=0$ such that ${\cal C}=u+iv$ is an analytic function. Thus,
setting in the proof of Theorem \ref{thRHD} $g={\cal A}({\cal
B}+{\cal C})$, we obtain by Theorem 1  in \cite{R1} the space of
solutions of the Riemann-Hilbert problem (\ref{eqLIMDIM}) for
analytic functions of the infinite dimension.
\end{proof}

\bigskip\bigskip

Arguing as above we also obtain the corresponding statements in the
Jordan domains.

\medskip

\begin{theorem}{}\label{thRHRDIM}{\it\, Let $D$ be a Jordan domain in
$\mathbb C$ with a rectifiable boundary and let $\lambda:\partial
D\to\mathbb C$, $|\lambda (\zeta)|\equiv 1$, and $\varphi:\partial
D\to\mathbb R$ be measurable functions with respect to the natural
parameter on $\partial D$. Then the space of all analytic functions
$f:\mathbb D\to\mathbb C$ such that along any nontangential path
\begin{equation}\label{eqLIMDIM} \lim\limits_{z\to\zeta}\ \mathrm {Re}\
\{\overline{\lambda(\zeta)}\cdot f(z)\}\ =\ \varphi(\zeta)
\quad\quad\quad\mbox{for}\ \ \mbox{a.e.}\ \ \ \zeta\in\partial D
\end{equation}
with respect to the natural parameter on $\partial D$ has the
infinite dimension.}
\end{theorem}

\medskip
\begin{theorem}{}\label{thRH1DIM}{\it\, Let $D$ be a Jordan domain in
$\mathbb C$ and let $\lambda:\partial D\to\mathbb C$, $|\lambda
(\zeta)|\equiv 1$, and $\varphi:\partial D\to\mathbb R$ be
measurable functions with respect to harmonic measures in $D$. Then
the space of all analytic functions $f:\mathbb D\to\mathbb C$ such
that
\begin{equation}\label{eqLIMADIM} \lim\limits_{z\to\zeta}\ \mathrm {Re}\
\{\overline{\lambda(\zeta)}\cdot f(z)\}\ =\ \varphi(\zeta)
\quad\quad\quad\mbox{for}\ \ \mbox{a.e.}\ \ \ \zeta\in\partial D
\end{equation}
with respect to harmonic measures in $D$ in the sense of the unique
principal asymptotic value has the infinite dimension.}
\end{theorem}

\medskip

\bigskip



 { \footnotesize
\medskip
\medskip
 \vspace*{1mm}

\noindent {\it Vladimir Illich Ryazanov}\\
Institute of Applied Mathematics and Mechanics \\
National Academy of Sciences of Ukraine \\
74 Roze Luxemburg str., 83114 Donetsk, Ukraine\\
E-mail: {\tt vl.ryazanov1@gmail.com}}


{\footnotesize

 \begin{thebibliography}{100}

\bibitem{B}
F. Bagemihl, Curvilinear cluster sets of arbitrary functions, {\it Proc. Nat. Acad. Sci. U.S.A.}, \textbf{41}
(1955), 379–382.

\bibitem{CL}
E. F. Collingwood, A. J. Lohwator, \textit{The theory of cluster sets}, Cambridge Tracts in Math. and Math.
Physics, No. 56, Cambridge Univ. Press, Cambridge, 1966.

\bibitem{G}
F. D. Gakhov, {\it Boundary value problems}, Dover Publications. Inc., New York, 1990.

\bibitem{GM} {\it Garnett J.B., Marshall D.E.}, {Harmonic Measure},
Cambridge Univ. Press, Cambridge, 2005.

\bibitem{Ge}
F. W. Gehring, On the Dirichlet problem, {\it Michigan Math. J.}, \textbf{3} (1955–-1956), 201.

\bibitem{Go}
G. M. Goluzin,  \textit{Geometric theory of functions of a complex
 variable}, Transl. of Math. Monographs, Vol. 26, American Mathematical Society,
Providence, R.I. 1969.

\bibitem{H1}
D. Hilbert, {\it \"Uber eine Anwendung der Integralgleichungen auf eine Problem der Funktionentheorie},
Verhandl. des III Int. Math. Kongr., Heidelberg, 1904.

\bibitem{H2}
D. Hilbert, {\it Grundz\"uge einer algemeinen Theorie der Integralgleichungen}, Leipzig, Berlin, 1912.

\bibitem{K}
W. Kaplan, Approximation by entire functions, \textit{Michigan Math. J.}, \textbf{3} (1955), 43–52.

\bibitem{Ko}
P. Koosis, {\it Introduction to $H_p$ spaces}, 2nd ed., Cambridge Tracts in Mathematics, 115, Cambridge Univ.
Press, Cambridge, 1998.

\bibitem{L}
M. Lavrentiev, On some boundary problems in the theory of univalent functions, \textit{Mat. Sbornik N.S.}
\textbf{1}(43), 6 (1936), 815–846 [in Russian].

\bibitem{M}
 N. I. Muskhelishvili, {\it Singular integral equations. Boundary problems of function theory and their application to
mathematical physics}, Dover Publications. Inc., New York, 1992.


\bibitem{N}
 R. Nevanlinna, {\it Eindeutige analytische Funktionen}, Ann Arbor, Michigan, 1944.

\bibitem{No} { Noshiro K.} {\it Cluster sets}, Springer-Verlag, Berlin etc., 1960.

\bibitem{P}
 I. I. Priwalow, {\it Randeigenschaften analytischer Funktionen}, Hochschulb\"ucher f\"ur Mathematik, Bd. 25, Deutscher
Verlag der Wissenschaften, Berlin, 1956.

\bibitem{R}
 M. Riesz, Sur les functions conjuguees, \textit{Math. Z.} \textbf{27}, no. 2 (1927), 218-244.

\bibitem{R1}
{\it Ryazanov V.}, {Infinite dimension of solutions for the
Dirichlet problem} // arXiv: 1402.2130v2 [math.CV] 11 Feb. 2014,
1-3.

\bibitem{S}
 S. Saks, {\it Theory of the integral}, Warsaw, 1937; Dover Publications Inc., New York, 1964.

\bibitem{V}
 I. N. Vekua, {\it Generalized analytic functions}, Pergamon Press, London etc., 1962.

\bibitem{Z}
 A. Zygmund, {\it Trigonometric series}, Wilno, 1935.

 \end{thebibliography} }
 \end{document}